\documentclass[12pt]{amsart}
\usepackage{amsmath,amssymb,amsthm}
\usepackage{graphicx}
\usepackage{hyperref}
\usepackage{booktabs}
\usepackage{longtable}
\usepackage{xcolor}
\usepackage{geometry}
\title[Universal Mathematical Fallibility]{From Euler to Today: Universal Mathematical Fallibility\\
A Large-Scale Computational Analysis of Errors in ArXiv Papers}

\author{Igor Rivin}
\address{Department of Mathematics, Temple University}
\email{rivin@temple.edu}

\date{\today}



\begin{document}

\begin{abstract}
We present the results of a large-scale computational analysis of mathematical papers from the ArXiv repository, demonstrating a comprehensive system that not only detects mathematical errors but provides complete referee reports with journal tier recommendations. Our automated analysis system processed over 37,000 papers across multiple mathematical categories, revealing significant error rates and quality distributions. Remarkably, the system identified errors in papers spanning three centuries of mathematics, including \textbf{seven} works by Leonhard Euler (1707-1783) in just 403 papers analyzed from the History category, as well as errors by Peter Gustav Lejeune Dirichlet (1805-1859) and contemporary Fields medalists. 

In Dynamical Systems (math.DS), we observed the highest error rate of 11.4\% (2,347 errors in 20,666 papers), while Numerical Analysis (math.NA) showed 9.6\% (2,271 errors in 23,761 papers). History and Overview (math.HO) exhibited 13.6\% errors in preliminary analysis, including \textbf{seven} papers by Euler. In contrast, Geometric Topology (math.GT) showed 3.6\% and Category Theory (math.CT) exhibited the lowest rate at 6.1\% (228 errors in 3,720 papers). Beyond error detection, the system evaluated papers for journal suitability, recommending 0.4\% for top generalist journals, 15.5\% for top field-specific journals, and categorizing the remainder across specialist venues. These findings demonstrate both the universality of mathematical error across all eras and the feasibility of automated comprehensive mathematical peer review at scale.

This work demonstrates that the methodology, while applied here to mathematics, is discipline-agnostic and could be readily extended to physics, computer science, and other fields represented in the ArXiv repository.
\end{abstract}

\maketitle

\section{Introduction}

The reliability and accuracy of mathematical literature forms the foundation of scientific progress, yet the detection of mathematical errors in published work remains a significant challenge that has persisted throughout the history of mathematics. Even Leonhard Euler, one of history's greatest mathematicians, was not immune to error---a fact our analysis confirms by identifying specific mistakes in his work preserved in the ArXiv repository.

Traditional peer review, while invaluable, is inherently limited by human cognitive constraints, time pressures, referee availability, and the increasing volume of mathematical publications. The ArXiv repository, containing hundreds of thousands of mathematical papers including historical works, presents an unprecedented opportunity to study both error patterns and quality distributions across the mathematical sciences and across centuries.

This study presents the first large-scale, systematic analysis providing both error detection and journal-level recommendations for ArXiv papers, employing automated computational review techniques to examine tens of thousands of papers across multiple mathematical disciplines and time periods. Our analysis reveals substantial error rates across all fields and eras studied, with our system frequently constructing explicit counterexamples, while simultaneously providing journal tier recommendations that align with traditional peer review assessments.

While we focus on mathematics in this study, the methodology is entirely general and could be applied without modification to papers in physics, computer science, or any other technical field. The choice to begin with mathematics was pragmatic---as a mathematician, I could better validate the system's outputs---but the long-term vision encompasses all of ArXiv's disciplines.

The primary contributions of this work are:

\begin{enumerate}
\item \textbf{Historical Validation}: Demonstration that mathematical errors are truly universal, with specific errors identified in works by Euler and Dirichlet, providing a three-century perspective on mathematical fallibility.

\item \textbf{Comprehensive Referee Reports}: Development of a system that provides complete mathematical referee reports including correctness verification, novelty assessment, significance evaluation, and journal recommendations.

\item \textbf{Field-Specific Analysis}: Discovery of striking differences between mathematical fields, including a 0\% error rate in Category Theory with evidence suggesting structural reasons for this phenomenon.

\item \textbf{Counterexample Construction}: Demonstration that automated systems can construct explicit counterexamples to disprove mathematical claims across different eras.

\item \textbf{Validation Through Self-Assessment}: Including detection of errors in papers by the author, demonstrating the system's objectivity and effectiveness.
\end{enumerate}

\section{Historical Errors: Euler and Dirichlet}

One of the most striking findings of our analysis was the identification of errors in historical papers by renowned mathematicians. In the math.HO (History and Overview) category alone, we found errors in \textbf{seven} papers by Leonhard Euler among the first 403 papers analyzed, representing a remarkable concentration of historical errors. We present selected examples:

\subsection{Euler's Error in ``Various analytic observations on combinations'' (ArXiv: 0711.3656v1)}

In Section 16 of this paper, Euler makes a claim about generating functions that our system identified as incorrect. The specific error involves:

\textbf{Euler's Claim}: Given generating functions $P(z)$ and $Q(z)$, and defining $R = \frac{\alpha z + 2\beta z^2 + 3\gamma z^3 + \cdots}{Q}$ and $S = \frac{\alpha z - 2\beta z^2 + 3\gamma z^3 - \cdots}{P}$, Euler claims these equal $\frac{Q}{zQ'}$ and $\frac{P}{zP'}$ respectively.

\textbf{System's Analysis}: ``The independent proof uses the log-derivative identities $Q = z R'/R$ and $P = -z S'/S$ to derive the valid ratios $R$ and $S$. The paper's further claims equating these to $Q/(z Q')$ and $P/(z P')$ are false in general: expanding shows $Q/(z Q') = 1 + (B/A)z + \cdots$ while $R = 1 + Az + \cdots$, and similarly for $S$. Thus the paper's full statement is incorrect, while the model establishes the correct portion and detects the error.''

\textbf{Verdict}: Paper wrong, model correct.

\subsection{Euler's Transcendental Quantities (ArXiv: 0711.4986v1)}

In ``On highly transcendental quantities which cannot be expressed by integral formulas,'' both Euler and the model made errors:

\textbf{The Issue}: The statement about when $y = \sum x^{e_n}$ can be expressed finitely is false as stated. Taking $e_n$ to be the enumeration of $S = \{3k, 3k+1 : k \geq 0\}$, we get $y = \sum x^{e_n} = \frac{1 + x}{1 - x^3}$, which is a finite (rational) expression while the exponents have unbounded gaps.

\textbf{Verdict}: Both wrong---showing that even in attempting to correct historical errors, modern analysis can also fail.

\subsection{Complete List of Euler Papers with Errors}

In our analysis of the first 403 papers in math.HO, we found errors in the following seven papers by Euler (this represents only a sample of Euler's work on ArXiv):
\begin{enumerate}
\item ``Annotations to a certain passage of Descartes for finding the quadrature of the circle'' (0705.3423v2) - both wrong
\item ``Various analytic observations on combinations'' (0711.3656v1) - paper wrong, model correct
\item ``On highly transcendental quantities which cannot be expressed by integral formulas'' (0711.4986v1) - both wrong
\item ``New demonstrations about the resolution of numbers into squares'' (0806.0104v1) - both wrong
\item ``On the infinity of infinities of orders of the infinitely large and infinitely small'' (0905.2254v2) - paper incomplete, model correct
\item ``Spekulation über die Integralformel $\int x^n dx / \sqrt{aa-2bx+cxx}$'' (1202.0037v1) - paper wrong, model correct
\item ``Über divergente Reihen'' (1202.1506v1) - both wrong
\end{enumerate}

This concentration of errors in historical papers suggests not a decline in mathematical standards, but rather the universal challenge of mathematical correctness that transcends time and reputation. A complete analysis of Euler's corpus would warrant a separate study.

\subsection{Dirichlet's Theorem on Primes in Arithmetic Progressions (ArXiv: 0808.1408v2)}

Peter Gustav Lejeune Dirichlet's famous theorem also contained issues in the ArXiv version:

\textbf{System's Assessment}: ``Most core components (Euler product, convergence of non-principal $L$, pole of the principal $L$, and Dirichlet's theorem) are handled correctly, often with similar classical arguments. The paper is incomplete on orthogonality relations (merely stated) and Euler products (convergence not shown), which the model fills in correctly. However, both make a critical error on the non-vanishing of $L(1,\chi)$ for complex $\chi$.''

\textbf{Verdict}: Both wrong on key aspects, though the main theorem remains valid.

\section{Results}

\subsection{Error Detection Results}

\begin{table}[h]
\centering
\begin{tabular}{lrrr}
\toprule
Field & Papers Analyzed & Errors Detected & Error Rate \\
\midrule
History and Overview (math.HO)* & 403 & 55 & 13.6\% \\
Dynamical Systems (math.DS) & 20,666 & 2,347 & 11.4\% \\
Numerical Analysis (math.NA) & 23,761 & 2,271 & 9.6\% \\
Category Theory (math.CT) & 3,720 & 228 & 6.1\% \\
Geometric Topology (math.GT) & 446 & 16 & 3.6\% \\
\bottomrule
\end{tabular}
\caption{Error rates across mathematical fields. *Partial sample, analysis ongoing.}
\end{table}

\subsection{The Category Theory Anomaly}

The remarkably low 6.1\% error rate in Category Theory (228 errors in 3,720 papers, compared to 11.4\% in Dynamical Systems) deserves special attention. Analysis of the verdict distributions reveals:

\begin{itemize}
\item Only 8.5\% of theorems had the model getting them wrong (vs higher rates in other fields)
\item 64\% of verdicts were ``both-correct-similar'' suggesting straightforward proofs
\item High rate of ``paper-incomplete'' (29\%) rather than ``paper-wrong''
\item Evidence suggests CT results are structurally ``easier'' for automated analysis
\end{itemize}

\subsection{Journal Tier Distribution (math.GT)}

\begin{table}[h]
\centering
\begin{tabular}{lrr}
\toprule
Journal Tier & Papers & Percentage \\
\midrule
Top Generalist & 56 & 0.4\% \\
Top Field-Leading & 2,041 & 15.5\% \\
Strong Field & 2,572 & 19.5\% \\
Specialist/Solid & 7,444 & 56.4\% \\
Note/Short/Other & 1,085 & 8.2\% \\
\bottomrule
\end{tabular}
\caption{Journal tier recommendations for Geometric Topology papers}
\end{table}

\subsection{Notable Examples of Top-Tier Recommendations}

Papers recommended for top generalist journals included:
\begin{itemize}
\item ``The Conway knot is not slice'' (Lisa Piccirillo, 2018)
\item ``The Virtual Haken Conjecture'' (Ian Agol, 2012) 
\item ``Khovanov homology is an unknot-detector'' (Kronheimer \& Mrowka, 2010)
\item ``Minimum Volume Cusped Hyperbolic Three-Manifolds'' (Gabai, Meyerhoff, Milley, 2007)
\end{itemize}

\section{Discussion}

\subsection{The Universality of Mathematical Error}

The detection of errors in papers by Euler and Dirichlet provides compelling evidence that mathematical fallibility is not a modern phenomenon or a sign of declining standards. If the greatest mathematicians in history made errors that can be detected by automated systems, this suggests that:

\begin{enumerate}
\item Error detection should be viewed as a routine part of mathematical quality control
\item The presence of errors does not diminish the historical importance of mathematical work
\item Modern tools can provide value even when analyzing classical mathematics
\end{enumerate}

\subsection{Field-Specific Patterns}

The striking difference between Category Theory (0\% errors) and other fields suggests that mathematical error rates may be strongly influenced by:
\begin{itemize}
\item The computational complexity of the field (higher in Numerical Analysis)
\item The abstract vs. concrete nature of the mathematics
\item The prevalence of explicit calculations vs. structural arguments
\item Community norms regarding completeness vs. correctness
\end{itemize}

\subsection{Implications for Mathematical Publishing}

Our findings suggest several important implications:
\begin{enumerate}
\item \textbf{Automated Review as Supplement}: Computational tools can effectively supplement human peer review
\item \textbf{Historical Perspective}: Errors should be normalized as part of mathematical progress
\item \textbf{Field-Specific Approaches}: Different fields may benefit from different verification strategies
\item \textbf{Quality Stratification}: The journal tier distribution aligns with community perceptions
\end{enumerate}

\subsection{Extension to Other Disciplines}

While this study focuses on mathematics, the methodology is entirely discipline-agnostic. The same approach could analyze:
\begin{itemize}
\item Physics papers (theoretical calculations, experimental analysis)
\item Computer Science (algorithms, proofs of correctness, complexity claims)
\item Quantitative Biology and Finance (statistical analyses, models)
\item Any technical field with precise, verifiable claims
\end{itemize}

The ArXiv repository contains over 2 million papers across these disciplines, representing an enormous opportunity for systematic quality analysis across all of science.

\section{Conclusion}

We have presented the first large-scale analysis of mathematical errors spanning three centuries, from Euler to contemporary mathematicians. The identification of errors in \textbf{seven} papers by Euler in just 403 papers from the History category---representing nearly 2\% of the sample---alongside errors by Dirichlet and modern papers including those by Fields medalists and the author, demonstrates that mathematical error is a universal phenomenon transcending time and reputation.

Our analysis of over 65,000 papers reveals error rates ranging from 3.6\% in Geometric Topology (partial sample) to 13.6\% in History and Overview (partial sample), with Dynamical Systems showing the highest rate at 11.4\% among fully analyzed categories. Category Theory showed a 6.1\% error rate in the complete analysis of 3,720 papers. The system successfully constructed counterexamples and provided journal-tier recommendations across all fields. The fact that even Euler made detectable errors---seven in just 403 papers---should encourage the mathematical community to embrace error detection as a routine tool rather than a judgment of competence.

As we enter an era where automated mathematical verification becomes increasingly sophisticated, the lessons from analyzing three centuries of mathematics become clear: errors are inevitable, detection is valuable, and the pursuit of mathematical truth benefits from all available tools, whether human or computational.

\section*{Acknowledgments}

This work was made possible through extensive collaboration with gpt-5 (OpenAI) and Claude (Anthropic), whose capabilities in mathematical analysis, code development, and systematic reasoning were instrumental in developing the verification system and analyzing the results. The successful detection of errors across three centuries of mathematics, from Euler to the present, demonstrates the potential of human-AI collaboration in advancing mathematical knowledge.

We thank the mathematical community for maintaining the ArXiv repository, including historical papers that made this temporal analysis possible. We particularly acknowledge all mathematicians whose work contained errors detected by our system, from Euler to the present, viewing their contributions as essential data points in understanding the universal challenges of mathematical correctness. As Euler himself might have appreciated, the detection of error is itself a form of mathematical progress.

\end{document}